\documentclass{article}

\usepackage{amsmath}
\usepackage{amssymb}
\usepackage{amsthm}
\usepackage{amscd}
\usepackage[all]{xy}

\newtheorem{theorem}{Theorem}[section]
\newtheorem{proposition}[theorem]{Proposition}
\newtheorem{lemma}[theorem]{Lemma}

\theoremstyle{definition}
\newtheorem{definition}[theorem]{Definition}

\theoremstyle{remark}

\newcommand{\Z}{\mathbb{Z}}
\newcommand{\HH}{\mathcal{H}}
\newcommand{\F}{\mathcal{F}}
\newcommand{\U}{\mathcal{U}}

\title{An approach toward a finite-dimensional definition of twisted $K$-theory}

\author{Kiyonori Gomi}

\date{}

\begin{document}

\maketitle

\begin{abstract}
This is an expository account of the following result: we can construct a group by means of twisted $\Z_2$-graded vectorial bundles which is isomorphic to  $K$-theory twisted by any degree three integral cohomology class.
\end{abstract}


\section{Introduction}

Topological $K$-theory admits a twisting by a degree three integral cohomology class. The resulting $K$-theory, known as \textit{twisted $K$-theory} \cite{A-S}, has its origin in the works of Donovan-Karoubi \cite{D-K} and Rosenberg \cite{R}, and has applications to D-brane charges \cite{B-M,K,W}, Verlinde algebras \cite{F-H-T} and so on.

\smallskip

As is well-known, ordinary $K$-theory has definitions by means of:
\begin{list}{}{\parsep=-2pt\topsep=4pt}
\item[(1)]
vector bundles;

\item[(2)]
the $C^*$-algebra of continuous functions; and

\item[(3)]
the space of Fredholm operators.
\end{list}
Twisted $K$-theory is usually defined by twisting the definitions (2) or (3). For a definition parallel to (1), there are the notions of \textit{twisted vector bundles} (\cite{L-U,M,T-X-U}, see also \cite{B-M,D-K,K}), and of \textit{bundle gerbe $K$-modules} \cite{BCMMS}. However, the definitions by means of these geometric objects are only valid for twisted $K$-theory whose ``twisting'', a third integral cohomology class, is of finite order: otherwise, there is no non-trivial such geometric object in finite-dimensions.

Toward a finite-dimensional definition of twisted $K$-theory valid for degree three integral cohomology classes of infinite order, we explain in this article Furuta's notion of generalized vector bundles \cite{F}, which we call \textit{vectorial bundles}. We also explain a notion of finite dimensional approximation of Fredholm operators, which provides a linear version of the finite dimensional approximation of the monopole equations \cite{F2}. We can use these notions to construct a group and an isomorphism from $K$-theory twisted by any degree three integral cohomology class. The proof of the result is only outlined. The detailed treatment will be provided elsewhere.

A possible application of the result above is to generalize the notions of \textit{2-vector bundles} \cite{B-D-R,B}. A 2-vector bundle of rank 1 due to Brylinski \cite{B} is a stack which reproduces the category of twisted vector bundles. Replacing the category of vector bundles by that of vectorial bundles, we can directly generalize the 2-vector bundles in \cite{B}. Similarly, we can also generalize the 2-vector bundles due to Baas, Dundas and Rognes \cite{B-D-R}, which they studied in an approach to geometric realization of elliptic cohomology. It seems interesting to apply their study of 2-vector bundles to the generalization of 2-vector bundles made of vectorial bundles.

\bigskip

\textit{Acknowledgments}.
I am grateful to M. Furuta for helpful discussions concerning this work. I am also indebted to T. Moriyama and A. Henriques for useful suggestions. I thank the organizers of School on Poisson geometry and related topics for the invitation to give a talk.


\section{Twisted $K$-theory and twisted vector bundles}

In this section, we review the definition of twisted $K$-theory by means of the space of Fredholm operators, following \cite{A-S}. We also review the notion of twisted vector bundles \cite{BCMMS,B-M,D-K,K,L-U,M,T-X-U}.

\smallskip

Unless otherwise mentioned, $X$ is a compact manifold through this article.


\subsection{Twisted $K$-theory}

The twisted $K$-theory we consider in this article is associated to a degree three integral cohomology class. To give the precise definition, we represent the class by a projective unitary bundle. Let $\HH$ be a separable Hilbert space of infinite dimension, and $PU(\HH)$ the projective unitary group $PU(\HH) = U(\HH)/U(1)$.

\begin{definition}
For a principal $PU(\HH)$-bundle over $X$, we define the \textit{twisted $K$-group} $K_P(X)$ to be the fiberwise homotopy classes of sections of the associated bundle $P \times_{Ad} \F(\HH)$ over $X$, where $PU(\HH)$ acts on the space $\F(\HH)$ of Fredholm operators on $\HH$ by adjoint.
\end{definition}

In the above definition, the topologies on $PU(\HH)$ and $\F(\HH)$ are understood to come from the operator norm. Notice that we can also use the compact-open topology in the sense of \cite{A-S}.

\medskip

If $P$ is a trivial bundle, then $K_P(X)$ is exactly the homotopy classes of continuous functions $X \to \F(\HH)$. Thus, in this case, we recover the $K$-group of $X$ by the well-known fact that $\F(\HH)$ is a classifying space of $K$-theory \cite{A}.

\smallskip

$PU(\HH)$-bundles over $X$ are classified by $H^3(X, \Z)$: since $U(\HH)$ is contractible by Kuiper's theorem, $PU(\HH)$ is homotopy equivalent to the Eilenberg-MacLane space $K(\Z, 2)$, so that the classifying space $BPU(\HH)$ is homotopy equivalent to $K(\Z, 3)$. If $P$ and $P'$ are isomorphic $PU(\HH)$-bundles, then the twisted $K$-groups $K_P(X)$ and $K_{P'}(X)$ are also (non-canonically) isomorphic. So we often speak of ``twisted $K$-theory twisted by a class in $H^3(X, \Z)$''.

\smallskip

We will call the cohomology class corresponding to $P$ the \textit{Dixmier-Douady class}, and denote it by $\delta(P) \in H^3(X, \Z)$. For later convenience, we recall here the construction of $\delta(P)$: take an open cover $\U = \{ U_\alpha \}$ of $X$ so that:
\begin{list}{$\bullet$}{\parsep=-2pt\topsep=4pt}
\item
there are local sections $s_\alpha : U_\alpha \to P|_{U_\alpha}$;

\item
there are lifts $g_{\alpha \beta} : U_{\alpha \beta} \to U(\HH)$ of the transition functions $\overline{g}_{\alpha \beta} : U_{\alpha \beta} \to PU(\HH)$.
\end{list}
Here we write $U_{\alpha \beta}$ for the overlap $U_\alpha \cap U_\beta$, and the transition function is defined by the relation $s_\alpha \overline{g}_{\alpha \beta} = s_\beta$.  Because of the cocycle condition for $\{ \overline{g}_{\alpha \beta} \}$, we can find a map $c_{\alpha \beta \gamma} : U_{\alpha \beta \gamma} \to U(1)$ such that $g_{\alpha \beta} g_{\beta \gamma} = c_{\alpha \beta \gamma} g_{\alpha \gamma}$. These maps comprise a \u{C}ech 2-cocycle $(c_{\alpha \beta \gamma}) \in \check{Z}^2(\U, \underline{U(1)})$ with its coefficients in the sheaf of germs of $U(1)$-valued functions, which represents $\delta(P)$ through the isomorphism $\check{H}^2(X, \underline{U(1)}) \cong H^3(X, \Z)$.


\subsection{Twisted vector bundles}

For a $PU(\HH)$-bundle $P$ over $X$, a \textit{twisted vector bundle} consists essentially of the data $(\U, E_\alpha, \phi_{\alpha \beta})$:
\begin{list}{$\bullet$}{\parsep=-2pt\topsep=4pt}
\item
an open cover $\U = \{ U_\alpha \}$ of $X$;

\item
vector bundles $E_\alpha$ over $U_\alpha$;

\item
isomorphisms of vector bundles $\phi_{\alpha \beta} : E_\beta|_{U_{\alpha \beta}} \to E_\alpha|_{U_{\alpha \beta}}$ over $U_{\alpha \beta}$ satisfying the ``twisted cocycle condition'' on $U_{\alpha \beta \gamma}$:
$$
\phi_{\alpha \beta} \phi_{\beta \gamma}
= c_{\alpha \beta \gamma} \phi_{\alpha \gamma},
$$
where $c_{\alpha \beta \gamma}$ is as in the previous subsection.
\end{list}
In the rigorous definition of twisted vector bundles, we have to include the choices of the local sections $s_\alpha$ and the lifts $g_{\alpha \beta}$. Though it is crucial to specify these choices in considering isomorphism classes of twisted vector bundles, we omit them for simplicity.

The isomorphism classes of twisted vector bundles $\mathrm{Vect}_P(X)$ gives rise to a semi-group by the direct sum of local vector bundles. Let $K(\mathrm{Vect}_P(X))$ denote the group given by applying the Grothendieck construction to $\mathrm{Vect}_P(X)$. Then the following fact is known. (See \cite{B-M,D-K,K,L-U,T-X-U}.)

\begin{proposition}
For a $PU(\HH)$-bundle whose Dixmier-Douady class $\delta(P)$ is of finite order, there exists an isomorphism:
$$
K_P(X) \longrightarrow K(\mathrm{Vect}_P(X)).
$$
\end{proposition}

Instead of twisted vector bundles, we can use bundle gerbe $K$-modules to obtain  an equivalent result \cite{BCMMS,C-W}.

\bigskip

The rank of a twisted vector bundle is a multiple of the order of $\delta(P)$. This can be seen readily as follows. Suppose that a twisted vector bundle has a finite rank $r$. Taking the determinant of the twisted cocycle condition, we have:
$$
\mathrm{det}\phi_{\alpha \beta} \mathrm{det}\phi_{\beta \gamma} =
c^r_{\alpha \beta \gamma} \mathrm{det}\phi_{\alpha \gamma}.
$$
Hence $(c^r_{\alpha \beta \gamma}) \in \check{Z}^2(\U, \underline{U(1)})$ is a coboundary and $r \delta(P) = 0$.

\smallskip

Because of the above remark, there are no no-trivial twisted vector bundles in the case where $\delta(P)$ is infinite order. So we cannot use twisted vector bundles of finite dimensions to realize $K_P(X)$ generally. In spite of this fact, collections of locally defined vector bundles seem to have the potential in defining $K_P(X)$ by means of finite-dimensional objects. An approach is to use the usual technique proving the isomorphism $K(X) \cong [X, \F(\HH)]$. In this approach, however, some complications prevent us from transparent management, in particular, in giving equivalence relation. The usage of Furuta's generalized vector bundle provides a more efficient approach, which we explain in the next section.


\section{Furuta's generalized vector bundle}

In this section, we explain a generalization of the notion of vector bundles introduced by M. Furuta \cite{F}. We call the generalized vector bundles \textit{vectorial bundles} for short. This notion is closely related to a finite-dimensional approximation of Fredholm operators. Applying these notions, we approach to our problem of defining twisted $K$-theory finite-dimensionally.


\subsection{Approximation of a Fredholm operator}
\label{subsec:single_approximation}

We begin with the simplest situation. A $\Z_2$-graded vectorial bundle over a single point is a pair $(E, h)$ consisting of:
\begin{list}{$\bullet$}{\parsep=-2pt\topsep=4pt}
\item
a $\Z_2$-graded Hermitian vector space $E = E^0 \oplus E^1$ of finite rank; and

\item
a Hermitian map $h : E \to E$ of degree 1.
\end{list}

By using a $\Z_2$-graded vectorial bundle over a point, we can approximate a single Fredholm operator as follows. Let $A : \HH \to \HH$ be a Fredholm operator. For simplicity, we assume the kernel or cokernel of $A$ is non-trivial We define the $\Z_2$-graded Hilbert space $\hat{\HH}$ by $\hat{\HH} = \HH \oplus \HH$, and the self-adjoint Fredholm operator $\hat{A} : \hat{\HH} \to \hat{\HH}$ of degree 1 by $\hat{A} = \left( \begin{array}{cc} 0 & A^* \\ A & 0 \end{array} \right)$. By the assumption, the spectrum $\sigma(\hat{A}^2)$ of the non-negative operator $\hat{A}^2$ contains 0. Since $\hat{A}^2$ is also Fredholm, $0 \in \sigma(\hat{A}^2)$ is a discrete spectrum. Hence there is a positive number $\mu$ such that:
\begin{list}{$\bullet$}{\parsep=-2pt\topsep=4pt}
\item
$\mu \not\in \sigma(\hat{A}^2)$;

\item
the subset $\sigma(\hat{A}^2) \cap [0, \mu)$ consists of a finite number of eigenvalues;

\item
the eigenspace $\mathrm{Ker}(\hat{A}^2 - \lambda)$ is finite dimensional for $\lambda \in \sigma(\hat{A}^2) \cap [0, \mu)$.
\end{list}

Let $0 = \lambda_1 < \lambda_2 < \cdots < \lambda_n < \mu$ be the distinct eigenvalues in $\sigma(\hat{A}^2) \cap [0, \mu)$. Then we have the following orthogonal decomposition of $\hat{\HH}$:
$$
\hat{\HH} =
(\hat{\HH}, \hat{A})_{\lambda_1} \oplus
(\hat{\HH}, \hat{A})_{\lambda_2} \oplus
\cdots \oplus
(\hat{\HH}, \hat{A})_{\lambda_n} \oplus
\hat{\HH}',
$$
where $(\hat{\HH}, \hat{A})_\lambda = \mathrm{Ker}(\hat{A}^2 - \lambda)$ is the eigenspace of $\hat{A}^2$ with its eigenvalue $\lambda$, and $\hat{\HH}'$ is the orthogonal complement. Notice that $\hat{A}$ preserves each eigenspace as well as the orthogonal complement. More precisely, $\hat{A}$ restricts to the trivial map on $(\hat{\HH}, \hat{A})_{\lambda_1} \cong \mathrm{Ker} A \oplus \mathrm{ker} A^*$, while $\hat{A}$ induces isomorphisms on $(\HH, \hat{A})_{\lambda_i}$ for $i > 1$ and $\hat{\HH}'$.

Now, cutting off the infinite-dimensional part $\hat{\HH}'$, we define $E$ and $h$ by $E = \bigoplus_i (\hat{\HH}, \hat{A})_{\lambda_i}$ and $h = \hat{A}|_E$. The pair $(E, h)$ is nothing but a $\Z_2$-graded vectorial bundle over a point.

\medskip

As a finite-dimensional approximation of a single Fredholm operator, we obtained a $\Z_2$-graded vectorial bundle over a single point. As will be explained in Subsection \ref{subsec:family_approximation}, a similar approximation is possible for a family of Fredholm operators parameterized by $X$. The resulting object is a $\Z_2$-graded vectorial bundle over $X$.


\subsection{Definition of vectorial bundle}

We now introduce $\Z_2$-graded vectorial bundles:

\begin{definition} \label{dfn:vectorial_bundle}
Let $(\U, (E_\alpha, h_\alpha), \phi_{\alpha \beta})$ be the following data:
\begin{list}{$\bullet$}{\parsep=-2pt\topsep=4pt}
\item
an open cover $\U = \{ U_\alpha \}$ of $X$;

\item
$\Z_2$-graded Hermitian vector bundles $E_\alpha$ over $U_\alpha$;

\item
Hermitian maps $h_\alpha : E_\alpha \to E_\alpha$ of degree 1;

\item
vector bundle maps $\phi_{\alpha \beta} : E_\beta|_{U_{\alpha \beta}} \to E_\alpha|_{U_{\alpha \beta}}$ of degree 0 over $U_{\alpha \beta}$ such that $h_\alpha \phi_{\alpha \beta} = \phi_{\alpha \beta} h_\beta$.
\end{list}
A \textit{$\Z_2$-graded vectorial bundle} over $X$ is defined to be data $(\U, (E_\alpha, h_\alpha), \phi_{\alpha \beta})$ satisfying the following conditions:
$$
\begin{array}{rcll}
\phi_{\alpha \alpha} &
\doteq & 1 &
\mbox{on $U_{\alpha}$}, \\
\phi_{\alpha \beta} \phi_{\beta \gamma} &
\doteq & \phi_{\alpha \gamma} &
\mbox{on $U_{\alpha \beta \gamma}$}. \\
\end{array}
$$
\end{definition}

In the above definition, the symbol \ $\doteq$ \ stands for an equivalence relation. The meaning of the first condition $\phi_{\alpha \alpha} \doteq 1$ is as follows.
\begin{quote}
For any point $x \in U_{\alpha}$, there are a neighborhood $V \subset U_{\alpha}$ of $x$ and a positive number $\mu$ such that: for all $y \in V$ and $v \in (E_\alpha, h_\alpha)_{y, < \mu}$ we have $\phi_{\alpha \alpha} (v) = v$.
\end{quote}
Here $(E_\alpha, h_\alpha)_{y, < \mu}$ is the subspace in the fiber of $E_\alpha$ at $y$ given by the direct sum of eigenspaces of $(h_\alpha)^2_y$ whose eigenvalues are less than $\mu$:
$$
(E_\alpha, h_\alpha)_{y, < \mu}
=
\bigoplus_{\lambda < \mu}
\mathrm{Ker} \left( (h_\alpha)^2_y - \lambda \right)
=
\bigoplus_{\lambda < \mu}
\left\{
v \in (E_\alpha)_y |\
(h_\alpha)_y^2 v = \lambda v
\right\}.
$$
The meaning of the second condition $\phi_{\alpha \beta} \phi_{\beta \gamma} \doteq \phi_{\alpha \gamma}$ is now obvious.

\medskip

\begin{definition}
Let $\mathbb{E} = (\U, (E_\alpha, h_\alpha), \phi_{\alpha \beta})$ and $\mathbb{E}' = (\U, (E'_\alpha, h'_\alpha), \phi'_{\alpha \beta})$ be $\Z_2$-graded vectorial bundles over $X$.

(a) A set $( f_\alpha )$ of vector bundle maps $f_\alpha : E_\alpha \to E'_\alpha$ of degree 0 such that $f_\alpha h_\alpha = h'_\alpha f_\alpha$ on $U_\alpha$ is said to be a \textit{homomorphism} from $\mathbb{E}$ to $\mathbb{E'}$, if we have $f_\alpha \phi_{\alpha \beta} \doteq \phi'_{\alpha \beta} f_\beta$ on $U_{\alpha \beta}$.

(b) A homomorphism $( f_\alpha ) : \mathbb{E} \to \mathbb{E}'$ is said to be an \textit{isomorphism}, if there exists a homomorphism $( f'_\alpha ) : \mathbb{E}' \to \mathbb{E}$ such that $f_\alpha f'_\alpha \doteq 1$ and $f'_\alpha f_\alpha \doteq 1$ on $U_\alpha$.
\end{definition}

In the above definition of homomorphism, $\mathbb{E}$ and $\mathbb{E'}$ share the same open cover. In the case where they have different open covers $\U$ and $\U'$ respectively, it suffices to take a common refinement of $\U$ and $\U'$.

\begin{definition}
A \textit{homotopy} between $\Z_2$-graded vectorial bundles $\mathbb{E}$ and $\mathbb{E}'$ over $X$ is defined to be a $\Z_2$-graded vectorial bundle $\tilde{\mathbb{E}}$ over $X \times [0, 1]$ such that $\mathbb{E}$ and $\mathbb{E}'$ are isomorphic to $\tilde{\mathbb{E}}|_{X \times \{ 0 \}}$ and $\tilde{\mathbb{E}}|_{X \times \{ 1 \}}$, respectively.
\end{definition}

We write $KF(X)$ for the set of homotopy classes of isomorphism classes of $\Z_2$-graded vectorial bundles. The set $KF(X)$ gives rise to a group by means of the direct sum of vector bundles given locally.

\smallskip

A $\Z_2$-graded (ordinary) vector bundle $E$ gives an example of a $\Z_2$-graded vectorial bundle by setting $\U = \{ X \}$ and $h = 0$. This construction induces a well-defined homomorphism $K(X) \to KF(X)$. In \cite{F}, Furuta proved:

\begin{proposition} \label{prop:Furuta}
The homomorphism $K(X) \to KF(X)$ is an isomorphism.
\end{proposition}


\subsection{Approximation of a family of Fredholm operators}
\label{subsec:family_approximation}

As a family version of the construction in Subsection \ref{subsec:single_approximation}, we can show:

\begin{lemma} \label{lem:approximation}
Let $A = \{ A_x \} : X \to \F(\HH)$ be a continuous map. For any point $p \in X$, there are a neighborhood $U_p$ of $p$ and a positive number $\mu_p$ such that the following family of vector spaces gives rise to a vector bundle over $U_p$:
$$
\bigcup_{x \in U_p} (\hat{\HH}, \hat{A}_x)_{< \mu_p}
=
\bigcup_{x \in U_p}
\bigoplus_{\lambda < \mu_p} \mathrm{Ker}(\hat{A}_x^2 - \lambda).
$$
\end{lemma}

A key to this lemma is that eigenvalues of $A_x$ is continuous in $x$.

\medskip

By means of the lemma, the family of Fredholm operators $A : X \to \F(\HH)$ yields a $\Z_2$-graded vectorial bundle $( \{ U_p \}_{p \in X}, (E_{U_p}, h_{U_p}), \phi_{U_p U_q})$, where the $\Z_2$-graded vector bundle $E_{U_p}$ is that in Lemma \ref{lem:approximation}, the Hermitian map $h_{U_p}$ is given by restricting the Fredholm operator: $h_{U_p}|_x = \hat{A}_x|_{E_{U_p}}$, and the map of vector bundles $\phi_{U_p U_q} : E_{U_q} \to E_{U_p}$ is the following composition of the natural inclusion and the orthogonal projection:
$$
\bigcup_{x \in U_p \cap U_q} (\hat{\HH}, \hat{A}_x)_{< \mu_q} \longrightarrow
(U_p \cap U_q) \times \hat{\HH} \longrightarrow
\bigcup_{x \in U_p \cap U_q} (\hat{\HH}, \hat{A}_x)_{< \mu_p}.
$$

\smallskip

The construction above induces a well-defined homomorphism
$$
\alpha : \  [X, \F(\HH)] \longrightarrow KF(X).
$$
This homomorphism is compatible with the isomorphism $\mathrm{ind} : [X, \F(\HH)] \to K(X)$ in \cite{A}. Namely, the following diagram is commutative:
$$
\begin{CD}
[X, \F(\HH)]  @= [X, \F(\HH)] \\
@V{\mathrm{ind}}VV @VV{\alpha}V \\
K(X) @>>> KF(X).
\end{CD}
$$
The compatibility follows from the fact that one can realize any vector bundle $E \to X$ as $E = \bigcup_{x \in X} \mathrm{Ker} \hat{A}_x$ by taking $A : X \to \F(\HH)$ such that $\sigma(\hat{A}^2_x) = \{ 0, 1 \}$. (See the proof of the surjectivity of $\mathrm{ind}$ in \cite{A}.)


\subsection{Twisted vectorial bundle and twisted $K$-theory}

We now apply vectorial bundles and finite dimensional approximations explained so far to twisted $K$-theory.

\smallskip

Recall that twisted vector bundles are defined by ``twisting'' the ordinary cocycle condition for vector bundles. In a similar way, for a $PU(\HH)$-bundle $P$, we define a \textit{twisted $\Z_2$-graded vectorial bundle} by replacing the ``cocycle condition'' $\phi_{\alpha \beta} \phi_{\beta \gamma} \doteq \phi_{\alpha \gamma}$ in Definition \ref{dfn:vectorial_bundle} by the ``twisted cocycle condition'':
$$
\phi_{\alpha \beta} \phi_{\beta \gamma}
\doteq
c_{\alpha \beta \gamma} \phi_{\alpha \gamma}.
$$

\medskip

A twisted $\Z_2$-graded vectorial bundle can be constructed from a section $\mathbb{A} : X \to P \times_{Ad} \F(\HH)$. The section gives a set of maps $\{ A_p : W_p \to \F(\HH) \}_{p \in X}$ such that $A_p = g_{pq} A_q g_{pq}^{-1}$, where $W_p$ is an open set containing $p$ and $g_{pq} : W_p \cap W_q \to U(\HH)$ is a lift of transition function of $P$. Now, we use Lemma \ref{lem:approximation} to define a Hermitian vector bundle over $U_p \subset W_p$ by $E_{U_p} = \bigcup_{x \in U_p} (\hat{\HH}, (\hat{A}_p)_x)_{< \mu_p}$. The map $A_p$ also defines a Hermitian map $h_{U_p}$ on $E_{U_p}$ by restriction. If we define $\phi_{U_p U_q} : E_{U_q} \to E_{U_p}$ by the following composition:
$$
\bigcup_{x \in U_{pq}} (\hat{\HH}, (\hat{A}_q)_x)_{< \mu_q} \to
U_{pq} \times \hat{\HH} \overset{\mathrm{id} \times g_{pq}}{\longrightarrow}
U_{pq} \times \hat{\HH} \to
\bigcup_{x \in U_{pq}} (\hat{\HH}, (\hat{A}_p)_x)_{< \mu_p},
$$
then $(\{ U_p \}, (E_{U_p}, h_{U_p}), \phi_{U_p U_q})$ is a twisted $\Z_2$-graded vectorial bundle.

\medskip

Introducing isomorphisms and homotopies in a similar way, we obtain the group $KF_P(X)$ of homotopy classes of isomorphism classes of twisted $\Z_2$-graded vectorial bundles. The above construction of twisted vectorial bundles induces the well-defined homomorphism
$$
\alpha : \ K_P(X) \longrightarrow KF_P(X).
$$

Since this map generalizes $\alpha : [X, \F(\HH)] \to KF(X)$, it is reasonable to expect that $\alpha$ gives rise to an isomorphism. In fact, we have:

\begin{theorem} \label{thm:main}
For any $PU(\HH)$-bundle $P$ over a compact manifold $X$, the homomorphism $\alpha : K_P(X) \longrightarrow KF_P(X)$ is bijective.
\end{theorem}

We sketch the proof of this result in the next subsection.


\subsection{Sketch of the proof of Theorem \ref{thm:main}}

The fundamental idea to prove Theorem \ref{thm:main} is to construct a kind of generalized cohomology theory on CW complexes by means of $KF_P(X)$.

\bigskip

As is known \cite{A-S,C-W}, the twisted $K$-group $K_P(X)$ fits into a certain generalized cohomology theory $\{ K^{-n}_P(X, Y) \}_{n \in \Z}$. In particular, for a CW pair $(X, Y)$ equipped with a $PU(\HH)$-bundle $P \to X$, we have the long exact sequence:
$$
\cdots \to
K^{-n-1}_{P|_Y}(Y) \overset{\delta_{-n-1}}{\to}
K^{-n}_P(X, Y) \to
K^{-n}_P(X) \to
K^{-n}_{P|_Y}(Y) \overset{\delta_{-n}}{\to} \cdots.
$$
Note that we can identify $K^{1}_P(X, Y)$ with $K^{-1}_P(X, Y) = K_{P \times I}(X \times I, Y \times I \cup X \times \partial I)$, and $\delta_0 : K^0_{P|_Y}(Y) \to K^1_P(X, Y)$ with the composition of the following maps:
$$
K^0_{P|_Y}(Y) \overset{\beta}{\longrightarrow}
K^0_{P|_Y \times D^2}(Y \times D^2, Y \times S^1) =
K^{-2}_{P|_Y}(Y) \overset{\delta_{-2}}{\longrightarrow}
K^{-1}_P(X, Y).
$$
Here $\beta$ induces the Bott periodicity, and is given by ``multiplying'' a map $T : D^2 \to \F(\HH)$ representing the generator of $K(D^2, S^1)$.

\bigskip

To construct a similar cohomology theory, we define $KF_P(X, Y)$ by using twisted $\Z_2$-graded vectorial bundles on $X$ whose support do not intersect $Y$. (The \textit{support} of a twisted $\Z_2$-graded vectorial bundle $\mathbb{E} = ( \mathcal{U}, (E_\alpha, h_\alpha), \phi_{\alpha \beta})$ on $X$ is the closure of the points $x \in X$ such that $(h_\alpha)_x$ is not invertible for an $\alpha$.) For $n \ge 0$, we put:
$$
KF^{-n}_P(X, Y) =
KF_{P \times I^n}(X \times I^n, Y \times I^n \cup X \times \partial I^n).
$$
Then $K^{-n}_P(X, Y)$ satisfies the (suitably modified) homotopy axiom and the excision axiom in the Eilenberg-Steenrod axioms. In a way parallel to the method in \cite{A}, we can also introduce a natural map $\delta_{-n} : KF^{-n}_{P|_Y}(Y) \to KF^{-n+1}_P(X, Y)$, and obtain the long exact sequence for a pair:
$$
\cdots \to
KF^{-1}_P(X) \to
KF^{-1}_{P|_Y}(Y) \overset{\delta_{-1}}{\to}
KF^0_P(X, Y) \to
KF^0_P(X) \to KF^0_{P|_Y}(Y).
$$
To extend this sequence, we put $KF^1_P(X, Y) = KF^{-1}_P(X, Y)$ and define $\delta_0 : KF^0_{P|_Y}(Y) \to KF^1_P(X, Y)$ to be the composition of:
$$
KF^0_{P|_Y}(Y) \overset{\beta}{\longrightarrow}
KF^0_{P|_Y \times D^2}(Y \times D^2, Y \times S^1) =
KF^{-2}_{P|_Y}(Y) \overset{\delta_{-2}}{\longrightarrow}
KF^{-1}_P(X, Y),
$$
where $\beta$ is given by tensoring a vector bundle representing the generator of $K(D^2, S^1)$. Then the composition of $KF^0_P(X) \to KF^0_{P|_Y}(Y) \to KF^1_P(X, Y)$ is trivial. (This sequence is not yet shown to be exact at this stage.)

\bigskip

The spaces $X \times I^n$ and $Y \times I^n \cup X \times \partial I^n$ in the definition of $KF^{-n}_P(X, Y)$ are used in that of $K^{-n}_P(X, Y)$. Hence the finite-dimensional approximation induces the natural homomorphism $\alpha_{-n} : K^{-n}_P(X, Y) \to KF^{-n}_P(X, Y)$ for $n \ge -1$. We can readily see that $\delta_{-n}$ ($n \ge 1$) commutes with $\alpha_{-n}$, since $\delta_{-n}$ is essentially defined by an inclusion map of spaces. If $X$ is compact, then $\beta$ commutes with $\alpha_{-n}$, so that $\delta_0$ does. The key to this fact is that the compactness allows us to choose a map $T : D^2 \to \F(\HH)$ realizing the generator of $K(D^2, S^1)$ in a way appropriate for the finite-dimensional approximation.

\bigskip

Now, for a finite CW complex $X$ and a $PU(\HH)$-bundle $P \to X$, we can prove the bijectivity of $\alpha_{-n} : K^{-n}_P(X) \to KF^{-n}_P(X)$ ($n \ge 0$) by the induction on the number of cells in $X$. Notice that, if $P$ is trivial, then an argument by using Proposition \ref{prop:Furuta} implies the bijectivity of $\alpha_{-n}$ ($n \ge 0)$. Thus, $\alpha_{-n}$ $(n \ge 0)$ is bijective in the case where $X$ is a point and has only one cell. If the number of the cells in $X$ is $r > 1$, then we can express $X$ as $X = Y \cup e^q$, where $Y \subset X$ is a subcomplex with $(r-1)$ cells, and $e^q \subset X$ is a cell of dimension $q$. We here make use of the commutative diagram ($n \ge 0$):
$$
\begin{array}{c@{}c@{}c@{}c@{}c@{}c@{}c@{}c@{}c}
K^{-n-1}_{P|_Y}(Y) & \ \to \ &
K^{-n}_P(X, Y) & \ \to \ &
K^{-n}_P(X) & \ \to \ &
K^{-n}_{P|_Y}(Y) & \ \to \ &
K^{-n+1}_P(X, Y) \\
\downarrow & &
\downarrow & &
\downarrow & &
\downarrow & &
\downarrow \\
KF^{-n-1}_{P|_Y}(Y) & \ \to \ &
KF^{-n}_P(X, Y) & \ \to \ &
KF^{-n}_P(X) & \ \to \ &
KF^{-n}_{P|_Y}(Y) & \ \to \ &
KF^{-n+1}_P(X, Y).
\end{array}
$$
We can assume that the first and the forth columns are bijective in the induction. The excision axiom implies that $K^{-n}_P(X, Y) \cong K^{-n}_P(D^q, S^{q-1})$ and $KF^{-n}_P(X, Y) \cong K^{-n}_P(D^q, S^{q-1})$. Since any $PU(\HH)$-bundle over $D^q$ is trivial, the second and fifth columns are also bijective. Thus, so is the third column. (The exactness at $KF^{-n}_{P|_Y}(Y)$ is not necessary in the five-term lemma.)


\end{document}